\documentclass[12pt,a4paper]{article}
\usepackage[cp1251]{inputenc} % cp866 - DOS, cp1251 - Windows
\usepackage{amsfonts}

\newcommand{\B}{$\hfill\Box$}
\newcommand{\al}{\alpha}

\newcommand{\ga}{\gamma}
\newcommand{\de}{\delta}
\newcommand{\la}{\lambda}
\newcommand{\om}{\omega}
\newcommand{\ee}{\varepsilon}
\newcommand{\vv}{\varphi}
\newcommand{\iy}{\infty}

\begin{document}

\vspace*{1.0cm}

\begin{center}
{\large\bf
SPECTRAL ANALYSIS FOR THE MATRIX STURM-LIOUVILLE OPERATOR ON A FINITE INTERVAL}
\\[0.2cm]
{\bf N. Bondarenko} \\[0.2cm]
\end{center}

\thispagestyle{empty}

{\bf Abstract.} The inverse spectral problem is investigated for the matrix Sturm-Liouville equation
on a finite interval. Properties of spectral characteristics are provided,
a constructive procedure for the solution of the inverse problem along with
necessary and sufficient conditions for its solvability is obtained.



{\bf 1. Introduction. }\\

{\it 1.1.}
In this paper, the inverse spectral problem is investigated
for the matrix Sturm-Liouville equation.
Inverse spectral problems are to recover operators from their spectral characteristics.

The scalar case has been studied fairly completely (see [1-3]).
The matrix case is more difficult for investigating.
Different statements of inverse spectral problems for the matrix case were given in [4], [5] and [6]
with corresponding uniqueness theorems.
A constructive solution procedure was provided in [7], but for the special case
of the simple spectrum only. Necessary and sufficient conditions were obtained in [8] for the case when
the spectrum is asymptotically simple, that is an important restriction.
Moreover, the method used by the authors of [8] does not give a reconstruction procedure.
We also note that necessary and sufficient conditions on spectral data were given in [9]
for Sturm-Liouville operators with matrix-valued potentials in the Sobolev space
$W_2^{-1}$. This class of potentials differs from one considered in this paper.

In this paper, we study the self-adjoint matrix Sturm-Liouville
operator in the general case, without any special restrictions on
the spectrum. Properties of spectral characteristics are
investigated, and necessary and sufficient conditions are obtained
for the solvability of the inverse problem. We provide a
constructive procedure for the solution of the inverse problem in
the general case, that is a generalization of the algorithm from
[7]. For solving the inverse problem we develop the ideas of the method
of spectral mappings [3].

\medskip

{\it 1.2.}
Consider the boundary value problem $L = L(Q(x), h, H)$ for the matrix Sturm-Liouville equation:
$$
    \ell Y: = -Y''+ Q(x) Y = \la Y, \quad x \in (0, \pi),                                           \eqno(1)
$$
$$
    U(Y) := Y'(0) - h Y(0) = 0, \quad V(Y) := Y'(\pi) + H Y(\pi) = 0.                               \eqno(2)
$$

Here $Y(x) = [y_k(x)]_{k = \overline{1, m}}$ is a column vector,
$\la$ is the spectral parameter, and $Q(x) = [Q_{jk}(x)]_{j, k =
\overline{1, m}}$, where $Q_{jk}(x) \in L_2(0, \pi)$ are
complex-valued functions. We will subsequently refer to the matrix
$Q(x)$ as the \textit{potential}. The boundary conditions are
given by the matrices $h = [h_{jk}]_{j, k = \overline{1, m}}$, $H
= [H_{jk}]_{j, k = \overline{1, m}}$, where $h_{jk}$ and $H_{jk}$
are complex numbers. In this paper we study the self-adjoint case,
when $Q = Q^*$, $h = h^*$, $H = H^*$.

Let $\vv(x, \la)$ and $S(x, \la)$ be matrix-solutions of equation
(1) under the initial conditions
$$
    \vv(0, \la) = I_m, \quad \vv'(0, \la) = h, \quad S(0, \la) = 0_m, \quad S'(0, \la) = I_m.
$$
where $I_m$ is the identity $m \times m$ matrix, $0_m$ is the zero $m \times m$ matrix.

The function $\Delta(\la) := \det[V(\vv)]$ is called the
\textit{characteristic function} of the boundary value problem
$L$. The zeros of the entire function $\Delta(\la)$ coincide with
the eigenvalues of $L$ (counting with their multiplicities, see
Lemma~3), and they are real in the self-adjoint case.

Let $\om = \om^*$ be some $m \times m$ matrix. We will write
$L(Q(x), h, H) \in A(\om)$, if the problem $L$ has a potential
from $L_2(0, \pi)$ and $h + H + \frac{1}{2} \int_0^{\pi} Q(x) \,
dx = \om$. Without loss of generality we may assume that $L \in
A(\om)$, where $\om \in D = \{\om \colon \om = \mbox{diag}
\{\om_1, \ldots, \om_m\}, \om_1 \le \ldots \le \om_m \}$.

One can achieve this condition applying the standard unitary transform.

In order to formulate the main result we need the following lemmas that will be proved in Section~2.

\medskip

{\bf Lemma 1.} { \it Let $L \in A(\om)$, $\om \in D$. The boundary
value problem $L$ has a countable set of eigenvalues $\{\la_{nq}
\}_{n \ge 0, q = \overline{1, m}}$, and
$$
     \rho_{nq} = \sqrt{\la_{nq}} = n + \frac{\om_q}{\pi n} + \frac{\kappa_{nq}}{n},
     \quad \{\kappa_{nq}\}_{n \ge 0} \in l_2, \, q = \overline{1, m}.                            \eqno(3)
$$
}

\medskip

Let $\Phi(x, \la) = [\Phi_{jk}(x, \la)]_{j, k = \overline{1, m}}$
be a matrix-solution of equation (1) under the boundary conditions
$U(\Phi) = I_m$, $V(\Phi) = 0_m$. We call $\Phi(x, \la)$ the
\textit{Weyl solution} for $L$. Put $M(\la) := \Phi(0, \la)$. The
matrix $M(\la) = [M_{jk}(\la)]_{j,k = \overline{1, m}}$ is called
the \textit{Weyl matrix} for $L$. The notion of the Weyl matrix
is a generalization of the notion of the Weyl function ($m$-function)
for the scalar case (see [1], [3]).
The Weyl functions and their generalizations often appear in applications and in pure
mathematical problems, and they are natural spectral characteristics in the inverse problem
theory for various classes of differential operators.

Using the definition for $M(\la)$ one can easily check that
$$
    M(\la) = -(V(\vv))^{-1} V(S).                                                           \eqno(4)
$$
It follows from representation (4) that the matrix-function
$M(\la)$ is meromorphic in $\la$ with simple poles in the
eigenvalues $\{\la_{nq}\}$ of $L$ (see Lemma~4).

Denote
$$
    \alpha_{nq} := \mathop{\mathrm{Res}}_{\la = \la_{nq}}M(\la).
$$

The data $\Lambda := \{\la_{nq}, \alpha_{nq} \}_{n \geq 0, \,q =
\overline{1, m}}$ are called the \textit{spectral data} of the
problem $L$.

Let $\{ \la_{n_k q_k} \}_{k \geq 0}$ be all the distinct
eigenvalues from the collection $\{ \la_{n q}\}_{n \geq 0, q =
\overline{1, m}}$. Put
$$
     \alpha'_{n_k q_k} := \alpha_{n_k q_k}, \, k \geq 0, \quad \alpha'_{nq} = 0_m,
     \, (n, q) \notin \{ (n_k, q_k )\}_{k \geq 0}.
$$

Denote
$$
1 = m_1 < \ldots < m_{p + 1} = m + 1,
$$
$$
\om_{m_s} = \ldots = \om_{m_{s + 1} - 1} =: \om^{(s)}, \quad s =
\overline{1, p}
$$
where $p$ is the number of different values among $\{\om_q\}_{q =
\overline{1, m}}$. Let
$$\alpha_n^{(s)} = \sum_{q = m_s}^{m_{s + 1} - 1} \alpha'_{nq}, \quad s = \overline{1, p}.$$

\medskip

{\bf Lemma 2.} {\it Let $L \in A(\om)$, $\om \in D$. Then the
following relation holds
$$
    \alpha_n^{(s)} = \frac{2}{\pi} I^{(s)} + \frac{\kappa_n^{(s)}}{n},
    \quad \{\| \kappa^{(s)}_{n} \|\}_{n \geq 0} \in l_2, \, s = \overline{1, p},                      \eqno(5)
$$
where
$$
    I^{(s)} = [I^{(s)}_{jk}]_{j, k = \overline{1, m}}, \quad I^{(s)}_{jk} = \left\{
                                                            \begin{array}{cc}
                                                               1, & m_s \leq j = k \leq m_{s+1}-1, \\
                                                               0, & otherwise,
                                                            \end{array} \right.
$$
and $\| . \|$ is a matrix norm: $\| a \| = \max_{j, k} a_{jk}$
}

\medskip

Consider the following inverse problem.

\medskip

{\bf Inverse Problem 1.} Given the spectral data $\Lambda$, construct $Q$, $h$ and $H$.

\medskip

We will write $\{\la_{nq}, \alpha_{nq} \}_{n \geq 0, q =
\overline{1, m}} \in \mbox{Sp}$, if for $\la_{nq} = \la_{kl}$ we
always have $\alpha_{nq} = \alpha_{kl}$.

The main result of this paper is

\medskip

{\bf Theorem 1.} {\it Let $\om \in D$. For data $\{\la_{nq},
\alpha_{nq}\}_{n \geq 0, q = \overline{1, m}} \in \mbox{Sp}$ to be
the spectral data for a certain problem $L \in A(\om)$ it is
necessary and sufficient to satisfy the following conditions.

1) The asymptotics (3) and (5) are valid.

2) All $\la_{nq}$ are real, $\alpha_{nq} = (\alpha_{nq})^*$,
$\alpha_{nq} \geq 0$ for all $n \geq 0$, $q = \overline{1, m}$,
and the ranks of the matrices $\alpha_{nq}$ coincide with the
multiplicities of $\la_{nq}$.

3) For any row vector $\gamma(\la)$ that is entire in $\la$, and
that satisfy the estimate
$$
    \gamma(\la) = O(\exp(|\mbox{Im}\, \sqrt{\la}|\pi)), \quad |\la| \to \iy,
$$
if $\gamma(\la_{nq})\alpha_{nq} = 0$ for all $n \geq 0$, $q =
\overline{1, m}$, then $\gamma(\la) \equiv 0$. }

\medskip

We prove necessity of the conditions of Theorem~1 in Section~2 and
sufficiency in
Section~4. In Section~3 the constructive procedure is provided for the solution of Inverse Problem~1.\\

{\bf 2. Necessity.}\\

{\it 2.1.}
Let us study some properties of the spectral data.

\medskip

{\bf Lemma 3.}
{\it
The zeroes of the characteristic function $\Delta(\la)$ coincide with the eigenvalues of the
boundary value problem $L$. The multiplicity of each zero $\la_0$ of the function
$\Delta(\la)$ equals to the multiplicity of the corresponding eigenvalue
(by the multiplicity of the eigenvalue we mean the number of the
corresponding linearly independent vector eigenfunctions).
}

\medskip

{\bf Proof.} 1) Let $\la_0$ be an eigenvalue of $L$, and let $Y^0$
be an eigenfunction corresponding to $\la_0$. Let us show that
$Y^0(x) = \vv(x, \la_0) Y^0(0)$. Clearly, $Y^0(0) = \vv(0, \la)
Y^0(0)$. It follows from $U(Y^0) = 0$ that ${Y^0}'(0) = h Y^0(0) =
\vv(0, \la) Y^0(0)$. Thus, $Y^0(x)$ and $\vv(x, \la_0)Y^0(0)$ are
the solutions for the same initial problem for the equation (1).
Consequently, they are equal.

2) Let us have exactly $k$ linearly independent eigenfunctions
$Y^1$, $Y^2$, \ldots, $Y^k$ corresponding to the eigenvalue
$\la_0$. Choose the invertible $m \times m$ matrix $C$ such that
the first $k$ columns of $\vv(x, \la_0) C$ coincide with the
eigenfunctions. Consider $Y(x, \la) := \vv(x, \la) C$, $Y(x, \la)
= [Y_q(x, \la)]_{q = \overline{1, m}}$, $Y_q(x, \la_0) = Y^q(x),
\, q = \overline{1, k}$. Clearly, that the zeros of $\Delta_1(\la)
:= \det V(Y) = \det V(\vv) \cdot \det C$ coincide with the zeros
of $\Delta(\la)$ counting with their multiplicities. Note that
$\la = \la_0$ is a zero of each of the columns $V(Y_1)$, \dots,
$V(Y_k)$. Hence, if $\la_0$ is the zero of the determinants
$\Delta_1(\la)$ and $\Delta(\la)$ with the multiplicity $p$, than $p
\geq k$.

3) Suppose that $p > k$. Rewrite $\Delta_1(\la)$ in the form
$$
    \Delta_1(\la) = (\la - \la_0)^k \Delta_2(\la),
$$
$$
    \Delta_2(\la) = \det \left[\frac{V(Y_1)}{\la-\la_0},
    \ldots, \frac{V(Y_k)}{\la-\la_0}, V(Y_{k + 1}), \ldots, V(Y_m)\right].
$$
In view our supposition, we have $\Delta_2(\la_0) = 0$, i.\,e. there exist
not all zero coefficients $\alpha_q$, $q = \overline{1, m}$
zero exist such that
$$
   \sum_{q = 1}^k \alpha_q \frac{d V(Y_q(x, \la_0))}{d\la} +
   \sum_{q = k + 1}^m \alpha_q V(Y_q(x, \la_0)) = 0.                                         \eqno(6)
$$

If $\alpha_q = 0$ for $q = \overline{1, k}$, then the function
$$
    Y^+(x, \la) := \sum_{q = k + 1}^{m} \alpha_q Y_q(x, \la)
$$
for $\la = \la_0$ is an eigenfunction corresponding to $\la_0$
that is linearly independent with $Y^q$, $q = \overline{1, k}$.
Since the eigenvalue $\la_0$ has exactly $k$ corresponding
eigenfunctions, we arrive at a contradiction.

Otherwise we consider the function
$$
    Y^+(x, \la) := \sum_{q = 1}^{k} \alpha_q Y_q(x, \la) +
    (\la - \la_0) \sum_{q = k + 1}^{m} \alpha_q Y_q(x, \la).
$$
It is easy to check that
$$
    \ell(Y^+(x, \la)) = \la Y^+(x, \la), \quad \ell\left(\frac{d}{d \la}Y^+(x, \la)\right) =
     \la \frac{d}{d \la}Y^+(x, \la) + Y^+(x, \la),
$$
$$
    U(Y^+) = U\left(\frac{d}{d \la} Y^+\right) = 0, \quad V(Y^+(x, \la_0)) = 0.
$$
Relation (6) is equivalent to the following one
$$
    V\left(\frac{d}{d \la} Y^+(x, \la_0)\right) = 0.
$$
Thus, we obtain that $Y^+(x, \la_0)$ is an eigenfunction, and
$\frac{d}{d \la}Y^+(x, \la_0)$ is a so-called \textit{associated
function} (see [10]) corresponding to $\la_0$. If we show that the
considered Sturm-Liouville operator does not have associated
functions, we will also arrive at a contradiction with
$\Delta_2(\la_0) \neq 0$, and finally, prove that $k = p$.

4) Let us prove that the self-adjoint operator given by (1), (2)
does not have associated functions. Let $\la_0$ be an eigenvalue
of $L$, and let $Y^0$ and $Y^1$ be a corresponding eigenfunction
and an associated function respectively, i.\,e. both $Y^0$ and
$Y^1$ satisfy (2) and
$$
    (\ell - \la_0) Y^0 = 0, \quad (\ell - \la_0) Y^1 = Y^0.
$$
This yields
$$
    ((\ell - \la_0)^2 Y^1, Y^1) = 0,
$$
for the scalar product defined by
$$
    (Y, Z) := \int_0^{\pi} Y^*(x) Z(x) \, dx.
$$
In case of the self-adjoint operator, we have $
    (\ell Y, Z) = (Y, \ell Z)
$
for any $Y$ and $Z$ satisfying (2), and the eigenvalue $\la_0$ is real.
Therefore,
$$
    ((\ell - \la_0) Y^1, (\ell - \la_0) Y^1) = (Y^0, Y^0) = 0,
$$
and $Y^0 = 0$. Recall that $Y^0$ is the eigenfunction, and get a contradiction.
\B

\medskip

{\bf Lemma 4.} {\it All poles of the Weyl matrix $M(\la)$ are
simple, and the ranks of the residue-matrices coincide with the
multiplicities of the corresponding eigenvalues of $L$. }

\medskip

{\bf Proof.} Let $\la_0$ be an eigenvalue of $L$ with a
multiplicity $k$, and let $Y_1$, $Y_2$, \ldots, $Y_k$ be linearly
independent vector eigenfunctions corresponding to $\la_0$.
Following the proof of Lemma~3, we introduce the invertible matrix
$C = [C_1, \ldots, C_m]$ such that $Y_q(x) = \vv(x, \la_0) C_q$,
$q = \overline{1, k}$. Consider the vector-function $Y(x, \la) =
\vv(x, \la) C$. Clearly, that $(V(\vv))^{-1} = C (V(Y))^{-1}$.
Write $V(Y(x, \la))$ in the form
\[
    V(Y(x, \la)) = [(\la - \la_0)W_1(\la), \ldots, (\la - \la_0) W_k(\la), W_{k + 1}(\la), \ldots, W_m(\la)],
\]
where
$$
      W_q(\la) = \frac{V(Y_q(x, \la))}{\la - \la_0}, \quad q = \overline{1, k},
$$
$$
      W_q(\la) = V(Y_q(\la)), \quad q = \overline{k + 1, m}.
$$
Clearly, $W_q(\la)$ are entire functions, and
$$\det W(\la) = \det [W_1(\la), \ldots, W_m(\la)] \neq 0$$
for $\la$ from a sufficiently small neighborhood of $\la_0$
(otherwise the multiplicity of the eigenvalue $\la_0$ is greater
than $k$). It is easy to show that
$$
    \det V(Y(x, \la)) = (\la - \la_0)^k \det W(\la),
$$
$$
     (V(Y(x, \la)))^{-1} = \left[\frac{X_1(\la)}{\la - \la_0}, \dots, \frac{X_k(\la)}{\la - \la_0}, X_{k+1}(\la), \dots, X_m(\la) \right]^t,
$$
where $X_q(\la)$ are analytic in a sufficiently small neighborhood
of $\la_0$ (the superscript $t$ stands for transposition). Using
(4) we get
$$
    \alpha_0 = \mathop{\mathrm{Res}}_{\la = \la_0} M(\la) = - \mathop{\mathrm{Res}}_{\la = \la_0} (V(\vv(x, \la)))^{-1} V(S(x, \la))
$$
$$
     = - \mathop{\mathrm{Res}}_{\la = \la_0} C \left[\frac{X_1(\la)}{\la - \la_0}, \dots, \frac{X_k(\la)}{\la - \la_0}, X_{k+1}(\la), \dots, X_m(\la) \right]^t
     V(S(x, \la))
$$
$$
     = - C \left[ X_1(\la_0), \dots, X_k(\la_0), 0, \dots, 0 \right]^t  V(S(x, \la_0))
    = - X V(S(x, \la_0)).
$$
Therefore, the poles of the Weyl matrix are simple, and $\mbox{\it rank} \; \alpha_0 \leq k$.

Let us prove the reverse inequality. Note that
\[
    \mathop{\mathrm{Res}}_{\la = \la_0} (V(\vv(x, \la)))^{-1} V(\vv(x, \la)) = 0_m = X V(\vv(x, \la_0)).
\]
Let $\psi(x, \la_0)$ be a solution of equation~(1) for $\la =
\la_0$ under the initial condition $V(\psi) = X^*$. Since columns
of the matrices $\vv(x, \la_0)$ and $S(x, \la_0) $ form a
fundamental system of solutions of equation (1), we have
\[
     \psi(x, \la_0) = \vv(x, \la_0) A + S(x, \la_0) B,
\]
\[
     X X^* = X V(\psi(x, \la_0)) = X V(\vv(x, \la_0)) A + X V(S(x, \la_0)) B = -\alpha_0 B.
\]
On the one hand, since $\det W \neq 0$, the vectors $X_q(\la_0)$
are linearly independent, therefore, $\mbox{\it rank} \; X X^* =
k$. On the other hand, $\mbox{\it rank} \; \alpha_0 B \leq
\mbox{\it rank} \; \alpha_0$. Thus, we conclude that $\mbox{\it
rank} \; \alpha_0 \geq k$. \B

\medskip

{\bf Lemma 5.} {\it Let $\la_0$, $\la_1$ be eigenvalues of $L$,
$\la_0 \neq \la_1$, and $\alpha_i =
\mathop{\mathrm{Res}}\limits_{\la = \la_{i}}M(\la)$, $i = 0, 1$.
The following relations hold
$$
    \alpha_0^* \int\limits_0^{\pi} \vv^*(x, \la_0) \vv(x, \la_0) \, dx \, \alpha_0 = \alpha^*_0,
$$
$$
    \alpha_0^* \int\limits_0^{\pi} \vv^*(x, \la_0) \vv(x, \la_1) \, dx \, \alpha_1 = 0_m.
$$
In particular, according to the first relation,
$$\alpha_0 = \alpha_0^* \geq 0.$$
}

{\bf Proof.}
Denote
$$
    \ell^*Z := -Z''+ Z Q(x), \quad
    V^*(Z) := Z'(\pi) + Z(\pi)H, \quad
    \langle Z, Y \rangle := Z'Y - Z Y',
$$
where $Z = [Z_k]^{t}_{k = \overline{1, m}}$ is a row vector ($t$ is the sign for the transposition).
Then $$\langle Z, Y \rangle_{x = \pi} =  V^*(Z) Y(\pi) - Z(\pi) V(Y).$$

If $Y(x, \la)$ and $Z(x, \mu)$ satisfy the equations $\ell Y(x, \la) =
\la Y(x,\la)$ and $\ell^*Z(x, \mu) = \mu Z(x, \mu)$, respectively,
then $\frac{d}{dx} \langle Z, Y \rangle = (\la - \mu)Z Y$. In
particular, if $\la = \mu$, then $\langle Z, Y \rangle$ does not
depend on $x$.

Since $\la_0$ is real, $\vv^*(x, \la_0)$ satisfies the equation
$\ell^* Z = \la_0 Z$. Hence,
$$
    \int_0^{\pi} \vv^*(x, \la_0) \vv(x, \la_0) \, dx =
    \lim_{\la \to \la_0} \frac{\langle \vv^*(x, \la_0), \vv(x, \la) \rangle|_0^{\pi}}{\la - \la_0}
$$
$$
    = \lim_{\la \to \la_0} \frac{V^*(\vv^*(x, \la_0)) \vv(x, \la) -
    \vv^*(x, \la_0) V(\vv(x, \la))}{\la - \la_0}.
$$
It follows from (4) and Lemma~4 that
$$
    V(\vv(x, \la_0)) \alpha_0 = - \lim_{\la \to \la_0} (\la - \la_0) V(\vv(x, \la))
    (V(\vv(x, \la)))^{-1} V(S(x, \la)) = 0_m.
$$
Analogously $\alpha_0^* V^*(\vv^*(x, \la_0)) = 0_m$. Consequently,
we calculate
$$
   \alpha_0^* \int\limits_0^{\pi} \vv^*(x, \la_0) \vv(x, \la_0) \, dx \, \alpha_0 =
   \alpha_0^* \vv^*(\pi, \la_0) \lim_{\la \to \la_0} \frac{V(\vv(x, \la))}{\la - \la_0}
$$
$$
    \times \lim_{\la \to \la_0} (\la - \la_0) (V(\vv(x, \la)))^{-1} V(S(x, \la))
   = \alpha_0^* \vv^*(\pi, \la_0) V(S(x, \la_0))
$$
$$
   = -\alpha_0^* \langle \vv^*(x, \la_0), S(x, \la_0) \rangle_{x = \pi}
   = -\alpha_0^* \langle \vv^*(x, \la_0), S(x, \la_0) \rangle_{x = 0} = \alpha_0^*.
$$
Similarly one can derive the second relation of the lemma.
\B

\medskip

{\it 2.2.} In this subsection we obtain asymptotics for the spectral data.

Denote $\rho := \sqrt \la$, $\mbox{Re} \, \rho \geq 0$, $\tau :=
\mbox{Im}\, \rho$, $G_{\de} = \{\rho \colon |\rho - k| \geq \de, k
= 0, 1, 2, \dots \}$, $\de > 0$. By the standard way (see
[3, Sec.~1.1]) one can obtain the estimate
$$
    \Delta(\la) = (- \rho \sin \rho \pi)^m + O(|\rho|^{m - 1}\exp(m |\tau| \pi))
    = (- \rho \sin \rho \pi)^m + (- \rho \sin \rho \pi)^{m - 1} O(\exp(|\tau|\pi))
$$
$$
    + \ldots + (-\rho \sin \rho \pi) O(\exp((m - 1)|\tau|\pi)) + O(\exp(m |\tau| \pi)),
    \quad |\rho| \to \iy.                                                                     \eqno(7)
$$

\medskip

{\bf Proof of Lemma~1.} 1) Consider the contour $\Gamma_N = \{\la
\colon |\la| = (N + 1/2)^2\}$. By virtue of~(7)
$$
   \Delta(\la) = f(\la) + g(\la), \quad
   f(\la) = (- \rho \sin \rho \pi)^m, \quad
   |g(\la)| \leq C |\rho|^{m - 1}\exp(m |\tau| \pi).
$$
If $\la \in \Gamma_N$ for sufficiently large $N$, we have
$|f(\la)| > |g(\la)|$. Then by Rouche's theorem the number of
zeros of $\Delta(\la)$ inside $\Gamma_N$ coincide with the number
of zeros of $f(\la)$ (counting with their multiplicities), i.\,e.
it equals $(N + 1) m$. Thus, in the circle $|\la| < (N + 1/2)^2$
there are exactly $(N + 1) m$ eigenvalues of $L$: $\{ \la_{nq}
\}_{n = \overline{0, N}, q = \overline{1, m}}$.

Applying Rouche's theorem to the circle $\gamma_n(\de) = \{\rho
\colon |\rho - n| \leq \de\}$, we conclude that for sufficiently
large $n$ there are exactly $m$ zeros of $\Delta(\rho^2)$ lying
inside $\gamma_n(\de)$, namely $\{\rho_{nq}\}_{q = \overline{1,
m}}$. Since $\de > 0$ is arbitrary, it follows that
$$
     \rho_{nq} = n + \ee_{nq}, \quad \ee_{nq} = o(1),
     \quad n \to \iy.                                                                         \eqno(8)
$$

Using (7) for $\rho = \rho_{nq}$, we get
$$
    (- \rho_{nq} \sin \rho_{nq} \pi)^m  + (- \rho_{nq} \sin \rho_{nq} \pi)^{m - 1} O(1) + \ldots +
    (-\rho_{nq} \sin \rho_{nq} \pi) O(1) + O(1) = 0, \quad n \to \iy.
$$
Denote $s_{nq} := |\rho_{nq} \sin \rho_{nq} \pi|$, and rewrite the obtained estimate in the form
$$
     s_{nq}^m \leq C_0 + C_1 s_{nq} + \ldots + C_{m - 1} s_{nq}^{m - 1}.                         \eqno(9)
$$
It follows from (9) that $s_{nq} \leq \max\{1, \sum_{k = 0}^{m - 1} C_k\}$. Otherwise we arrive at a
contradiction:
$$
   s_{nq}^m > \sum_{k = 0}^{m - 1} C_k s_{nq}^{m - 1} \geq \sum_{k = 0}^{m - 1} C_k s_{nq}^k.
$$
Hence, $|\rho_{nq} \sin \rho_{nq} \pi| \leq C$. Using (8) we get
$$
    \sin \rho_{nq} \pi = \sin \ee_{nq} \pi \cos n \pi = O(n^{-1}), \quad \ee_{nq} = O(n^{-1}), \quad n \to \iy.
$$
Together with (8) this yields
$$
    \rho_{nq} = n + O(n^{-1}), \quad n \to \iy.
$$

2) Let us derive the more precise asymptotic formula.
One can easily show that
$$
     V(\vv) = - \rho \sin \rho \pi \cdot I_m + \om \cos \rho \pi + \kappa(\rho),
$$
where
$$
     \kappa(\rho) = \frac{1}{2} \int\limits_{0}^{\pi}Q(t) \cos \rho(\pi - 2 t)\,dt +
     O\left(\frac{1}{\rho}\exp(|\tau|\pi)\right).
$$

Consider the linear mappings $z_n(\rho)$ that map the circles
$\{\rho \colon |\rho - n| \leq C / n\}$
(note that $\rho_{nq}$ lie in these circles for a fixed sufficiently large $C$)
to the circle $\{z \colon |z| \leq R\}$:
$$
    \rho = n + \frac{z_n(\rho)}{\pi n}.
$$
For $|z| \leq R$ we have
$$
     V(\vv) = (-1)^n (\om - z_n(\rho)I_m + \kappa_n(z_n(\rho))).                       \eqno(10)
$$
Using the representation for $\kappa(\rho)$ we get $\kappa_n(z) =
o(1)$, $n \to \iy$, uniformly with respect to $z$ in the circle
$\{ z \colon |z| \leq R\}$. Moreover, for each sequence
$\{z^0_n\}_{n\geq 0} \subset \{z \colon |z| \leq R \}$ we have
$\{ \| \kappa_n(z^0_n) \| \}_{n \geq 0} \in l_2$ and $\sum\limits_{n \geq
0} \|\kappa_n(z^0_n)\|^2 < C$, where $C$ is some constant.
Consequently,
$$\Delta(\rho^2) = \pm f(z_n(\rho)) + g_n(z_n(\rho)),$$
where $f(z) = \det (\om - z I_m)$, $g_n(z) = o(1)$, $n \to \iy$
(uniformly with respect to $z \in \{z \colon |z| \leq R\}$), and
the choice of sign $\pm$ depends only on $n$. Fix $0 < \de < 1/2
\min\limits_{q, l \colon \om_q \neq \om_l}|\om_q - \om_l|$ and
introduce the contours $\gamma_q = \{z \colon |z - \om_{q}| =
\de\}$. Clearly, the inequality $|f(z)| > |g_n(z)|$ holds on
$\gamma_q$ for all sufficiently large $n$, and by Rouche's theorem
two analytic functions $\Delta(\rho_n^2(z))$ and $f(z)$ have an
equal number of zeros inside $\gamma_q$ (here $\rho_n$ is the
inverse mapping to $z_n$). Thus, we have
$$
      \rho_{nq} = n + \frac{\om_q}{\pi n} + \frac{\kappa_{nq}}{n},
      \quad \kappa_{nq} = o(1), \, n \to \iy, \quad q = \overline{1, m}.
$$
Substituting this formula into (10) we get
\[
      V(\vv) = (-1)^n (\om - \om_q I_m - \pi \kappa_{nq} I_m + \kappa_n(z_n(\rho_{nq}))).
\]
Since $\{ \| \kappa_n(z_n(\rho_{nq})) \|\} \in l_2$, one can easily prove that $\{ \kappa_{nq} \} \in l_2$.
\B

\medskip

{\bf Proof of Lemma~2.} 1) Let $\tilde M(\la)$ be the Weyl matrix
for the problem $\tilde L(\tilde Q, \tilde h, \tilde H)$, such
that $\tilde Q(x) = \frac{2}{\pi}\om$, $\tilde h = \tilde H =
0_m$. Then $\tilde \alpha^{(s)} = \frac{2}{\pi} I^{(s)}$, $s =
\overline{1, p}$.

Consider the contours $\gamma^{(s)}_{n} = \{\la \colon |\la - (n^2
+ \frac{2}{\pi}\om^{(s)})| = R\}$, $R = \frac{1}{\pi}
\min\limits_{q, l \colon \om_q \neq \om_l}|\om_q - \om_l|$. Using
the residue theorem and taking Lemma~1 into account, we deduce
$$
   \frac{1}{2 \pi i} \int\limits_{\gamma_n^{(s)}}(M(\la) - \tilde M(\la))\,d \la =
   \sum\limits_{q = m_s}^{m_{s+1}-1}\alpha'_{nq} - \sum\limits_{q = m_{s}}^{m_{s+1}-1}\tilde \alpha'_{nq} =
    \alpha_n^{(s)} - \frac{2}{\pi}I^{(s)},
   n \geq n^*, \, s = \overline{1, p}.
$$
One can easily show that $M_{jk}(\la) = -\frac{\Delta_{jk}(\la)}{
\Delta(\la)}$, where
$$
    \Delta_{jk}(\la) = \det[V(\vv_1), \dots, V(\vv_{j - 1}),
    V(S_k), V(\vv_{j + 1}), \dots, V(\vv_m)].
$$
Using this representation, we arrive at
$$
     M_{jk}(\la) - \tilde M_{jk}(\la) =
     \frac{\Delta(\la) \tilde \Delta_{jk}(\la)
     - \Delta_{jk}(\la)\tilde \Delta(\la)}{\Delta(\la)
     \tilde \Delta(\la)},
     \quad j, k = \overline{1, m}.                                                            \eqno(11)
$$

Let us use the mappings $z_n$ introduced in the proof of Lemma~1:
\[
     \rho = n + \frac{z_n(\rho)}{\pi n}.
\]
If $\la \in \gamma_n^{(s)}$, then $0 < \de_1 \leq|z_n(\rho) -
\om_q|$ for all $q = \overline{1, m}$, and $|z_n(\rho) -
\om^{(s)}| \leq \de_2$. Hence, the estimate for $\Delta(\la)$
obtained in the proof of Lemma~1 is valid: $\Delta(\la) = \pm
f(z_n(\rho)) + o(1)$, $\la \in \gamma_n^{(s)}$, $n \to \iy$
(uniformly with respect to $\la$).

Similarly, we estimate
$$
    \Delta_{jk}(\la) = \pm \frac{f(z_n(\rho))}{z_n(\rho) - \om_j} + o(1) \quad for\, j = k,
$$
$$
    \Delta_{jk}(\la) = o(1) \quad for \,j \neq k,
$$
$$
\la \in \gamma_n^{(s)}, \, n \to \iy, \quad j, k = \overline{1,
m}.
$$
Convergence of the remainders is uniform with respect to $\la$,
the choice of sign $\pm$ depends only on $n$. Analogous relations
hold for $\tilde \Delta(\la)$ and $\tilde \Delta_{jk}(\la)$.

Substituting these estimates into~(11) and taking into account
that $C_1 \leq |f(z_n(\rho))| \leq C_2$ for $\la \in
\gamma_n^{(s)}$, we arrive at
$$
    M_{jk}(\la) - \tilde M_{jk}(\la) = o(1), \quad j, k = \overline{1, m}, \quad \la \in \gamma_n^{(s)},
$$
$$
   \frac{1}{2 \pi i} \int\limits_{\gamma_n^{(s)}} (M(\la) - \tilde M(\la))\, d\la = o(1), \quad
            \alpha_n^{(s)} = \frac{2}{\pi} I^{(s)} + \eta_n^{(s)},
    \, \eta_n^{(s)} = o(1), \quad n \to \iy.
$$

2) Below one and the same symbol $\{\kappa_n\}$ denotes various matrix sequences such that $ \{ \| \kappa_n\|\} \in l_2$.
Using the standard asymptotics
$$
     \vv(x, \la) = \cos \rho x \cdot I_m + Q_1(x) \frac{\sin \rho x}{\rho} + \int\limits_0^x \frac{\sin \rho (x - 2 t)}{2 \rho }Q(t) \, dt
     + O\left(\frac{\exp |\tau|x}{\rho}\right),
$$
$$
    |\rho| \to \iy, \quad x \in [0, \pi],
$$
where $Q_1(x) = h + \int\limits_0^x Q(t) \, dt$, one can easily show that
$$
     \int\limits_0^{\pi} \vv^*(x, \la_{nq}) \vv(x, \la_{nl}) \, dx = \frac{\pi}{2}I_m + \frac{\kappa_n}{n}, \quad \la_{nq} - \la_{nl} = \frac{\kappa_n}{n}.
$$
Applying Lemma~5, we get
\[
     \alpha_{nq} \left(\frac{\pi}{2}I_m + \frac{\kappa_n}{n}\right) \alpha_{nq} = \alpha_{nq}, \quad n \geq 0, \quad q = \overline{1, m}.
\]
Clearly, $\|\alpha_{nq} \| \leq C$, $n \geq 0$, $q = \overline{1, m}$.
Consequently, $\frac{\pi}{2} \alpha^2_{nq} = \alpha_{nq} + \frac{\kappa_n}{n}$.
Similarly we derive $\alpha_{nq} \alpha_{nl} = \frac{\kappa_n}{n}$, $m_s \leq q, l \leq m_{s+1}-1$, $q \neq l$,
$s = \overline{1, p}$.
Thus,
$$
     \frac{\pi}{2}(\alpha^{(s)}_{n})^2 =\frac{\pi}{2}
     \left(\sum\limits_{q = m_s}^{m_{s+1}-1} \alpha'_{nq} \right)^2 = \frac{\pi}{2}
     \sum\limits_{q = m_s}^{m_{s+1}-1}(\alpha'_{nq})^2 + \frac{\kappa_n}{n}
    = \sum\limits_{q = m_s}^{m_{s+1}-1} \alpha'_{nq} + \frac{\kappa_n}{n} =
     \alpha^{(s)}_n + \frac{\kappa_n}{n}.
$$
Substitute the result of point~1 into this equality:
$$
   \frac{\pi}{2}\left(\frac{2}{\pi} I^{(s)} + \eta_n^{(s)}\right)^2 = \frac{2}{\pi} I^{(s)} + \eta_n^{(s)} + \frac{\kappa_n}{n},
$$
$$
   (I_m - 2 I^{(s)}) \eta_n^{(s)} = \frac{\pi}{2}(\eta_n^{(s)})^2 + \frac{\kappa_n}{n}.
$$
Consequently, $\eta_n^{(s)} = \frac{\kappa_n}{n}$.
\B

\medskip

{\it 2.3.} {\bf Proof of Theorem~1 (necessity).} The first two
conditions are fulfilled by Lemmas~1, 2, 4, 5.

Let $\gamma(\la)$ be a function described in condition~3. Recall
that
$$V(\vv(x, \la_{nq})) \alpha_{nq} = 0_m.$$
Since
$$\mbox{\it rank}\; V(\vv(x, \la_{nq})) + \mbox{\it rank}\; \alpha_{nq} = m$$
 and $\gamma(\la_{nq})\alpha_{nq} = 0$, we get
$\gamma(\la_{nq}) = C_{nq} V(\vv(x, \la_{nq}))$, i.\,e. the row
$\gamma(\la_{nq})$ is a linear combination of the rows of the
matrix $V(\vv(x, \la_{nq})$ (here $C_{nq}$ is a row of
coefficients). Consider
$$f(\la) = \gamma(\la) (V(\vv(x, \la)))^{-1}.$$
The matrix-function $(V(\vv(x, \la)))^{-1}$ has simple poles in
$\la = \la_{nq}$, therefore, we calculate
$$
   \mathop{\mathrm{Res}}_{\la = \la_{nq}} f(\la) = \gamma(\la_{nq}) \mathop{\mathrm{Res}}_{\la = \la_{nq}}(V(\vv(x, \la)))^{-1}
$$
$$
   = C_{nq} \lim_{\la \to \la_{nq}} V(\vv(x, \la)) \lim_{\la \to \la_{nq}}(\la - \la_{nq})
  (V(\vv(x, \la)))^{-1} = 0.
$$
Hence, $f(\la)$ is entire. It is easy to show that
\[
    \| (V(\vv(x, \la)))^{-1}\| \leq C_{\de} |\rho|^{-1} \exp(-|\tau| \pi), \quad \rho \in G_{\de},
\]
where $G_{\de} = \{\rho \colon |\rho - k| \geq \de, k = 0, 1,
2, \dots \}$, $\de > 0$. From this we conclude that $\|
f(\la)\| \leq \frac{C}{|\rho|}$ in $G_{\de}$. By the maximum
principle this estimate is valid in the whole $\la$-plane. Using
Liouville`s theorem, we obtain $f(\la) \equiv 0$. Consequently,
$\gamma(\la) \equiv 0$. \B

\medskip

Note that in the scalar case condition~3 follows from the first two conditions of Theorem~1.
Indeed, in the scalar case, we have $\ga(\la_n) \al_n = 0$, $n \geq 0$, where
$\al_n$ are positive real numbers. Hence, $\ga(\la_n) = 0$.
Having the spectrum $\{ \la_n\}_{n \geq 0}$ we
can construct the characteristic function (see [3, Theorem~1.1.4]):
$$
    \Delta(\la) = \pi (\la - \la_0) \prod_{n = 1}^{\iy} \frac{\la_n - \la}{n^2},
$$
and using asymptotics (3) for the eigenvalues we get the estimate
$$
    \| \Delta(\la)\| \geq C_{\de} |\rho| \exp(|\tau| \pi), \quad \rho \in G_{\de}.
$$
Then we introduce $f(\la) = \frac{\gamma(\la)}{\Delta(\la)}$ and follow the proof of
necessity in Theorem~1.

In the general case, condition~3 is essential and cannot be omitted, that is shown by the following example.

\medskip

{\bf Example 1.} Let $m = 2$, $\la_{01} \neq \la_{02}$, $\la_{n1}
= \la_{n2} =n^2$, $n \geq 1$,
$$
    \al_{01} = \alpha_{02} = \left[ \begin{array}{cc} \frac{1}{\pi} & 0 \\ 0 & 0 \end{array} \right], \quad
    \al_{n1} = \alpha_{n2} = \left[ \begin{array}{cc} \frac{2}{\pi} & 0 \\ 0 & \frac{2}{\pi} \end{array} \right], \, n \geq 1.
$$
The data $\{\la_{nq}, \al_{nq} \}$ satisfy conditions~1-2 of
Theorem~1. Let us show that they do not satisfy condition~3, and
consequently, they cannot be spectral data of $L$. The relations
$\ga(\la_{nq})\alpha_{nq} = 0$, $n \geq 0$, $q = \overline{1,
m}$ for this example can be rewritten in the form $\ga(\la) =
[\ga_1(\la), \ga_2(\la)]$, $\ga_1(\la_{01}) =
\ga_1(\la_{02}) = \ga_1(n^2) = 0$, $\ga_2(n^2) = 0$, $n
\geq 1$. Clearly, if we put $\ga_1(\la) = 0$, $\ga_2(\la) =
\frac{\sin \rho \pi}{\rho}$,
we arrive at a contradiction with condition~3.

\medskip

Below we investigate condition~3 in some special cases.

\medskip

{\bf Example 2 (full multiplicities).} Let $\la_{n1} = \la_{n2} = \ldots = \la_{nm} =: \la_n$ for all
$n \geq 0$. Then $\mbox{\it rank} \; \al_{nq} = m$, and each of the linear systems
$\ga(\la_{nq}) \al_{nq} = 0$ has the unique solution $\ga(\la_n) = 0$. We get the situation similar
to the scalar case, because in view of asymptotics (3), $\{ \la_n\}_{n \geq 0}$ can
be treated as eigenvalues of some scalar problem. Therefore, condition~3 holds automatically.

\medskip

We will say that the relations $\ga(\la_{nq}) \al_{nq} = 0$, $q = \overline{1, m}$ are {\it separated} for
some fixed $n$, if they yield $\ga_q(\la_{nq}) = 0$ for all $q = \overline{1, m}$. For example,
they are separated in the case of full multiplicities, or when
the matrices $\al_{nq}$ have a proper diagonal form.

\medskip

{\bf Example 3.} Let the relations $\ga(\la_{nq}) \al_{nq} = 0$ be separated for all $n > n_0$.
Then each component $\ga_q(\la)$ has zeros $\{\la_{nq}\}_{n > n_0}$.
If $\ga(\la)$ is the function from condition~3, each $\ga_q(\la)$
cannot have more than $n_0$ additional zeros (counting with their multiplicities).
Otherwise we consider its zeros as the eigenvalues of a scalar problem and prove
that $\ga_q(\la) \equiv 0$.

If $\ga(\la)$ is entire, and $\gamma(\la) = O(\exp(|\mbox{Im}\, \sqrt{\la}|\pi)), \quad |\la| \to \iy$,
its order is not greater than $1/2$. Therefore, by Hadamard`s factorization theorem
$\ga_q(\la)$ can be presented in the form
$$
    \ga_q(\la) = (C_{q0} + C_{q1} \la + C_{q2} \la^2 + \ldots + C_{q, n_0} \la^{n_0}) P_q(\la),
    \quad P_q(\la) = \prod_{n > n_0} \left(1 - \frac{\la}{\la_{nq}} \right).
$$
We substitute this factorization into  $\ga(\la_{nq}) \al_{nq} = 0$, $n \leq n_0$, $q = \overline{1, m}$,
and obtain the system of linear equations with respect to $C_{q0}$, $C_{q1}$, \dots, $C_{q n_0}$, $q = \overline{1, m}$.

More precisely, let $\la_1$, \dots, $\la_N$ be the first $N = (n_0 + 1) m$ eigenvalues, and let
$\al_1$, \dots, $\al_N$ be the corresponding residue-matrices.
For each $j = \overline{1, N}$, we choose a non-zero column $v_j$ of $\al_j$.
In case of a group of multiple values among $\la_j$, $j = \overline{1, N}$, they have
a common matrix $\alpha_j$, and its rank equals their multiplicity, and we choose
linearly independent columns. Consider $N \times N$ matrix $P$ with the columns
$$[v_{j1} P_1(\la_j), v_{j1} \la_j P_1(\la_j), \dots, v_{j1} \la_j^{n_0} P_1(\la_j),
\ldots, v_{jm} P_m(\la_j), v_{jm} \la_j P_m(\la_j), \dots, v_{jm} \la_j^{n_0} P_m(\la_j)],$$
$j = \overline{1, N}$. Clearly, that the condition $\ga(\la_{nq}) \al_{nq} = 0$,
$n \leq n_0$, $q = \overline{1, m}$ is equivalent to the linear system with
the matrix $P$. Each solution of this system corresponds to $\ga(\la)$,
satisfying condition~3 of Theorem~1. Thus, the condition~3 is fulfilled
iff the determinant of $P$ is not zero. \\

{\bf 3. Solution of Inverse Problem~1.}\\

{\it 3.1.}
Let the spectral data $\Lambda$ of the boundary value problem $L
\in A(\om)$, $\om \in D$, be given.

Denote
$$
     D(x, \la, \mu) = \frac{ \langle \vv^*(x, \bar{\mu}),
     \vv(x, \la)\rangle}{\la - \mu} = \int\limits_0^x
     \vv^*(t, \bar{\mu}) \vv(x, \la) \, dt.                                       \eqno(12)
$$

We choose an arbitrary model boundary value problem $\tilde L =
L(\tilde Q(x), \tilde h, \tilde H) \in A(\om)$ (for example, one
can take $\tilde Q(x) = \frac{2}{\pi} \om$, $\tilde h = 0_m$,
$\tilde H = 0_m$). We agree that if a certain symbol $\gamma$
denotes an object related to $L$, then the corresponding symbol
$\tilde \gamma$ with tilde denotes the analogous object related to
$\tilde L$. Put
$$
   \xi_n = \sum_{q = 1}^m |\rho_{nq} - \tilde \rho_{nq}| + \sum_{s = 1}^p \sum_{q = m_s}^{m_{s+1}-1} |\rho_{nq} - \rho_{n m_s}| +
    \sum_{s = 1}^p \sum_{q = m_s}^{m_{s+1}-1} |\tilde \rho_{nq} - \tilde \rho_{n m_s}| + \sum_{s = 1}^p \| \alpha_n^{(s)} - \tilde \alpha_n^{(s)}\|.
$$
According to Lemmas~1 and 2,
\[
    \Omega := \left( \sum_{n = 0}^{\iy} ((n + 1) \xi_n)^2\right)^{1/2} < \iy, \quad \sum_{n = 0}^{\iy} \xi_n < \iy.
\]
Denote
$$
    \begin{array}{c}
    \la_{nq0} = \la_{nq}, \quad \la_{nq1} = \tilde \la_{nq}, \quad \rho_{nq0} = \rho_{nq}, \quad \rho_{nq1} = \tilde \rho_{nq}, \quad
    \alpha'_{nq0} = \alpha'_{nq}, \quad \alpha'_{nq1} = \tilde \alpha'_{nq}, \\
    \vv_{nqi}(x) = \vv(x, \la_{nqi}), \quad
    \tilde \vv_{nqi}(x) = \tilde \vv(x, \la_{nqi}), \\
    F_{klj, nqi}(x) =
    \alpha'_{klj} D(x, \la_{nqi}, \la_{klj}), \quad
    \tilde F_{klj, nqi}(x) = \alpha'_{klj}\tilde D(x, \la_{nqi}, \la_{klj}), \\
    n, k \geq 0, \quad q, l = \overline{1, m}, \quad i, j = 0, 1.
    \end{array}
$$

By the standard way (see [3, Lemma~1.6.2]), using Schwarz's lemma, we get

\medskip

{\bf Lemma 6.}
{\it
The following estimates are valid for $x \in [0, \pi]$, $n, k \geq 0$, $r, s = \overline{1, m}$, $m_r < q < m_{r+1}$, $m_s < l < m_{s+1}$, $i, j = 0, 1$:
$$
    \begin{array}{c}
   \|\vv_{nqi}(x)\| \leq C, \quad
   \|\vv_{n m_r 0}(x) - \vv_{n m_r 1}(x) \| \leq C \xi_n, \\
   \|\vv_{nqi}(x) - \vv_{nm_ri}(x) \| \leq C \xi_n, \quad
   \|F_{klj, nqi}(x) \| \leq \frac{C}{|n - k| + 1}, \\
   \left\|\sum\limits_{l = m_s}^{m_{s+1}-1} (F_{kl0, nm_r1}(x) - F_{kl1, nm_r1}(x)) \right\| \leq \frac{C \xi_k}{|n - k| + 1}, \\
   \|F_{klj, nqi}(x) - F_{klj, nm_ri}(x) \|, \, \|F_{klj, nm_r0}(x) - F_{klj, nm_r1}(x) \| \leq  \frac{C \xi_n}{|n - k| + 1}, \\
   \Biggl\| \sum\limits_{l = m_s}^{m_{s+1}-1} (F_{kl0, nqi}(x) - F_{kl0, nm_ri}(x) - F_{kl1, nqi}(x) + F_{kl1, nm_ri}(x)) \Biggr\|
   \leq \frac{C \xi_n \xi_k} {|n - k| + 1}, \\
   \Biggl\|\sum\limits_{l = m_s}^{m_{s+1}-1} (F_{kl0, nm_r0}(x) - F_{kl0, nm_r1}(x) - F_{kl1, nm_r0}(x) + F_{kl1, nm_r1}(x)) \Biggr\| \leq \frac{C \xi_n \xi_k} {|n - k| + 1}.
   \end{array}
$$
The analogous estimates are also valid for $\tilde \vv_{nqi}(x)$,
$\tilde F_{klj, nqi}(x)$. }

\medskip

The lemma similar to the following one has been proved in [7] by the contour integral method.

\medskip

{\bf Lemma 7.}
{\it
The following relations hold
$$
    \tilde \vv(x, \la) = \vv(x, \la) + \sum_{k = 0}^{\iy}
    \sum_{l = 1}^m (\vv_{kl0}(x) \al'_{kl0} \tilde D(x, \la, \la_{kl0}) -
    \vv_{kl1}(x) \al'_{kl1} \tilde D(x, \la, \la_{kl1}))                                  \eqno(13)
$$
$$
    \tilde D(x, \la, \mu) - D(x, \la,\mu) = \sum_{k = 0}^{\iy}
    \sum_{l = 1}^m (D(x, \la_{kl0}, \mu) \tilde D(x, \la, \la_{kl0}) -
    D(x, \la_{kl1}, \mu) \tilde D(x, \la, \la_{kl1}).
$$
Both series converge absolutely and uniformly with respect to $x \in [0, \pi]$ and
$\la$, $\mu$ on compact sets.
}

\medskip

Analogously one can obtain the following relation
$$
    \tilde \Phi(x, \la) = \Phi(x, \la) +
    \sum_{k = 0}^{\iy} \sum_{l = 1}^m \sum_{j = 0}^1
    (-1)^j \vv_{klj}(x) \al'_{klj} \frac{\langle
    \tilde \vv^*_{klj}(x), \tilde \Phi(x, \la)
         \rangle}{\la - \la_{klj}}.                                                           \eqno(14)
$$

It follows from Lemma~7 that
$$
    \tilde \vv_{nqi}(x) = \vv_{nqi}(x) + \sum_{k = 0}^{\iy}
    \sum_{l = 1}^m (\vv_{kl0} \tilde F_{kl0, nqi}(x) -
    \vv_{kl1} \tilde F_{kl1, nqi}(x)),                                                    \eqno(15)
$$
$$
    \tilde F_{\eta p \om, nqi}(x) - F_{\eta p \om, nqi}(x) =
    \sum_{k = 0}^{\iy} \sum_{l = 1}^m (F_{\eta p \om, kl0}(x) \tilde F_{kl0, nqi}(x)
    -  F_{\eta p \om, kl1}(x) \tilde F_{kl1, nqi}(x))                                     \eqno(16)
$$
for $n, \eta \geq 0$, $q, p = \overline{1, m}$, $i, \om  = 0, 1$.

Denote
$$
     \ee_0(x) = \sum_{(k, l, j) \in V} (-1)^j \vv_{klj}(x) \alpha'_{klj} \tilde \vv^*_{klj}(x),
     \quad \ee(x) = -2 \ee_0'(x).                                                        \eqno(17)
$$
Using (5) and Lemma~6 one can easily check that the series in (17)
converges absolutely and uniformly on $[0, \pi]$, and the
function $\ee_0(x)$ is absolutely continuous, and the components of
$\ee(x)$ belong to $L_2(0, \pi)$.

\medskip

{\bf Lemma 8.}
{\it
The following relations hold
$$
    Q(x) = \tilde Q(x) + \ee(x), \quad
    h = \tilde h - \ee_0(0), \quad H = \tilde H + \ee_0(\pi),                           \eqno(18)
$$
}

{\bf Proof.}
Differentiating (13) twice with respect to $x$ and using (12) and (17) we get
$$
    \tilde \vv'(x, \la) - \ee_0(x) \tilde \vv(x, \la) =
    \vv'(x, \la) + \sum_{k = 0}^{\iy} \sum_{l = 1}^m \sum_{j = 0}^1
    (-1)^j \vv'_{klj}(x) \al'_{klj} \tilde D(x, \la, \la_{klj}),
$$
$$
    \tilde \vv''(x, \la) = \vv''(x, \la) + \sum_{k = 0}^{\iy}
    \sum_{l = 1}^m \sum_{j = 0}^1
    (-1)^j [\vv''_{klj}(x) \al'_{klj} \tilde D(x, \la, \la_{klj})
$$
$$
    + 2 \vv'_{klj}(x) \al'_{klj} \tilde \vv^*_{klj}(x) \tilde \vv(x, \la) +
    \vv_{klj}(x) \al'_{klj} (\tilde \vv^*_{klj}(x) \tilde \vv(x, \la))'].
$$
We replace here the second derivatives, using equation (1), and then replace $\vv(x, \la)$,
using (13). This yields
$$
    \tilde Q(x) \vv(x, \la) = Q(x) \tilde \vv(x, \la) +
    \sum_{k = 0}^{\iy} \sum_{l = 1}^m \sum_{j = 0}^1 (-1)^j
    [\vv_{klj}(x) \al'_{klj} \langle \tilde \vv^*_{klj}(x), \tilde \vv(x, \la)\rangle
$$
$$
    + 2 \vv'_{klj}(x) \al'_{klj} \tilde \vv^*_{klj}(x) \tilde \vv(x, \la) +
    \vv_{klj}(x) \al'_{klj} (\tilde \vv^*_{klj}(x) \tilde \vv(x, \la))'].
$$
Cancelling terms with $\tilde \vv'(x, \la)$ we arrive at $Q(x) =
\tilde Q(x) + \ee(x)$.

Further,
$$
    \tilde \vv'(0, \la) - (h + \ee_0(0)) \tilde \vv(0) =
    U(\vv) + \sum_{k = 0}^{\iy} \sum_{l = 1}^m \sum_{j = 0}^1
    (-1)^j U(\vv_{klj}) \al'_{klj} D(0, \la, \la_{klj}) = 0_m.
$$
Since $\tilde \vv(0, \la) = I_m$, $\tilde \vv'(0, \la) = \tilde h$, we
obtain $h = \tilde h - \ee_0(0)$.

Similarly, using (14) one can get
$$
    \tilde \Phi'(\pi, \la) + (H - \ee_0(\pi)) \Phi(\pi, \la) =
    V(\Phi) + \sum_{k = 0}^{\iy} \sum_{l = 1}^m \sum_{j = 0}^1
    (-1)^j V(\vv_{klj}) \al'_{klj} \frac{\langle \tilde \vv^*_{klj}(x), \tilde \Phi(x, \la) \rangle_{|x = \pi}}
    {\la - \la_{klj}}.
$$
For $j = 0$ we have $V(\vv_{kl0}) \al'_{kl0} = 0_m$.
For $j = 1$
$$
    \langle \tilde \vv^*_{kl1}(x), \tilde \Phi(x, \la) \rangle_{|x = \pi} =
    \tilde V^*(\tilde \vv^*_{kl1}) \tilde \Phi(\pi, \la) -
    \tilde \vv^*_{kl1}(\pi) \tilde V(\tilde \Phi).
$$
Recall that $V(\Phi) = 0_m$, $\tilde V(\tilde \Phi) = 0_m$ and
$\al'_{kl1} \tilde V^*(\tilde \vv^*_{kl1}) = 0_m$. Consequently,
we arrive at $\tilde \Phi'(\pi, \la) + (H - \ee_0(\pi)) \Phi(\pi, \la) = 0_m$.
Together with $\tilde V(\tilde \Phi) = 0_m$ this yields
$H = \tilde H + \ee(\pi)$.
\B

\medskip

For each fixed $x \in [0, \pi]$, the relation (15) can be
considered as a system of linear equations with respect to
$\vv_{nqi}(x)$, $n \geq 0$, $q = \overline{1, m}$, $i = 0, 1$. But
the series in (15) converges only ``with brackets''. Therefore, it
is not convenient to use (15) as a main equation of the inverse
problem. Below we will transfer (15) to a linear equation in a
corresponding Banach space of sequences.\\

{\it 3.2.}
Denote $\chi_n := \xi_n^{-1}$ for $\xi_n \neq 0$ and $\chi_n = 0$ for $\xi_n = 0$.
Let $V$ be a set of indices $u = (n, q, i)$, $n \geq 0$, $q = \overline{1, m}$, $i = 0, 1$. For
each fixed $x \in [0, \pi]$, we define the row-vector $\psi(x) = [\psi_u(x)]_{u \in V}$
and the matrix $R(x) = [R_{v, u}(x)]_{v, u \in V}$,
$v = (k, l, j)$, $u = (n, q, i)$, by the formulae
$$
    \left.
    \begin{array}{c}
    \psi_{n m_s 0}(x) = \chi_n(\vv_{n m_s 0}(x) - \vv_{n m_s 1}(x)), \quad \psi_{n m_s 1}(x) = \vv_{n m_s 1}(x), \\
    \psi_{nqi}(x) = \chi_n(\vv_{nqi}(x) - \vv_{n m_s i}(x)), \\
    R_{k m_s 0, n m_r 0}(x) = \chi_n \xi_k \sum\limits_{l = m_s}^{m_{s+1}-1} (F_{kl0, n m_r 0}(x) - F_{kl0, n m_r 1}(x)), \\
    R_{k m_s 0, n m_r 1}(x) = \xi_k \sum\limits_{l = m_s}^{m_{s+1}-1} F_{kl0, n m_r 1}(x), \\
    R_{k m_s 0, nqi}(x) = \chi_n \xi_k \sum\limits_{l = m_s}^{m_{s+1}-1} (F_{kl0, nqi}(x) - F_{kl0, n m_r i}(x)), \\
    R_{klj, n m_r 0}(x) = (-1)^j \chi_n \xi_k (F_{klj, n m_r 0}(x) - F_{klj, n m_r 1}(x)), \\
    R_{klj, n m_r 1}(x) = (-1)^j \xi_k F_{klj,  m_r 1}(x), \\
    R_{klj, nqi}(x) = (-1)^j \chi_n \xi_k (F_{klj, nqi}(x) - F_{klj, n m_r i}(x)), \\
    R_{k m_s 1, n m_r 0}(x) = \chi_n  \sum\limits_{l = m_s}^{m_{s+1}-1} (F_{kl0, nm_r0}(x) - F_{kl0, nm_r1}(x) \\ - F_{kl1, nm_r0}(x) + F_{kl1, nm_r1}(x)), \\
    R_{k m_s 1, nqi}(x) = \chi_n \sum\limits_{l = m_s}^{m_{s+1}-1} (F_{kl0, nqi}(x) - F_{kl0, nm_ri}(x) - F_{kl1, nqi}(x) + F_{kl1, nm_ri}(x)), \\
    R_{k m_s 1, n m_r 1}(x) = \sum\limits_{l = m_s}^{m_{s+1}-1} (F_{kl0, nm_r1}(x) - F_{kl1, nm_r1}(x)), \\
    n, k \geq 0, \quad r, s = \overline{1, p}, \quad m_s < l < m_{s + 1}, \quad m_r < q < m_{r + 1}.
    \end{array}
    \right\}                                                                      \eqno(19)
$$
Analogously we define $\tilde \psi(x)$, $\tilde R(x)$ by replacing
in the previous definitions $\vv_{nqi}(x)$ by $\tilde
\vv_{nqi}(x)$ and $F_{klj, nqi}(x)$ by $\tilde F_{klj, nqi}(x)$.

We will also use a shorter notation. Consider the row vectors with
matrix components
$$
    \varphi_n(x) = [\varphi_{n10}(x), \varphi_{n11}(x), \varphi_{n20}(x), \varphi_{n21}(x), \dots, \varphi_{nm0}(x), \varphi_{nm1}(x)],
$$
$$
    \psi_n(x) = [\psi_{n10}(x), \psi_{n11}(x), \psi_{n20}(x), \psi_{n21}(x), \dots, \psi_{nm0}(x), \psi_{nm1}(x)], \quad n \geq 0,
$$
and defined analogously $2 m \times 2 m$ matrices
$F^{-}_{k,n}(x)$, $R_{k,n}(x)$, $n, k \geq 0$,
$F^{-}_{klj, nqi}(x) = (-1)^j F_{klj, nqi}(x)$.
Then definitions (19) of $\psi_{nqi}(x)$ and $R_{klj, nqi}(x)$ can be rewritten in the form
$$
   \psi_n = \varphi_n X_n, \quad R_{k,n} = X^{-1}_k F^{-}_{k,n} X_n, \quad n, k \geq 0.      \eqno(20)
$$
where $X_n$ is a $2 m \times 2 m$ matrix with components determined from (19).
Analogously we define $\tilde \varphi_n(x)$, $\tilde \psi_n(x)$ and $\tilde F^{-}_{k, n}(x)$, $\tilde R_{k, n}(x)$.
Now we can rewrite (15) and (16) in the form
$$
    \tilde \varphi_n = \varphi_n + \sum_{k = 0}^{\iy}
    \varphi_k \tilde F^{-}_{k, n}, \quad n \geq 0,                                            \eqno(21)
$$
$$
    \tilde F^-_{\eta, n} - F^-_{\eta, n} = \sum_{k = 0}^{\iy} F^-_{\eta,k} \tilde F^-_{k, n}  \eqno(22)
$$

By virtue of Lemma~6
$$
   \|\psi_{nqi}(x) \|, \, \|\tilde \psi_{nqi}(x) \| \leq C,
$$
$$
   \|R_{klj, nqi}(x) \|, \, \|\tilde R_{klj, nqi}(x) \| \leq \frac{C \xi_k}{|n - k| + 1},     \eqno(23)
$$
where $C$ does not depend on $x, n, q, i, k, l, j$

Let $a_u$, $u \in V$, be $m \times m$ matrices. Consider the
Banach space $B$ of bounded sequences $a = [a_u]_{u \in V}$ with
the norm $\|a \|_B = \sup\limits_{u \in V} \|a_u \|$. It follows
from (23) that for each fixed  $x \in [0, \pi]$, the operators $I +
\tilde R(x)$ and $I - R$ (here I is the identity operator), acting from $B$ to
$B$, are linear bounded operators.

\medskip

{\bf Theorem 2.}
{\it
For each fixed $x \in [0, \pi]$, the vector $\psi(x) \in B$
satisfies the equation
$$
    \tilde \psi (x) = \psi (x) (I + \tilde R(x))                                              \eqno(24)
$$
in Banach space $B$. Moreover, the operator $I + \tilde R(x) $
has a bounded inverse operator, i.\,e. equation (24) is uniquely solvable.
}

\medskip

{\bf Proof.}
Using (20) we get
$$
    \varphi_n = \psi_n X^{-1}_n, \quad F^{-}_{k, n} = X_k R_{k, n} X^{-1}_n,
$$
Substituting these relations into (21), we derive
$$
     \tilde \psi_n X^{-1}_n = \psi_n X^{-1}_n + \sum_{k = 0}^{\infty} \psi_k X^{-1}_k X_k \tilde R_{k, n} X^{-1}_n =
     \psi_n X^{-1}_n + \sum_{k = 0}^{\infty} \psi_k \tilde R_{k, n} X^{-1}_n \quad n \geq 0.
$$
Multiplying the result by $X_n$, we arrive at (24).

Similarly we get from (22) that
$$
    \tilde R_{\eta,n} - R_{\eta, n} = \sum_{k = 0}^{\iy} R_{\eta,k} \tilde R_{k,n}.
$$
This yields $\tilde R(x) - R(x) - R(x) \tilde R(x) = 0$, i.\,e.
$(I - R(x))(I + \tilde R(x)) = I$. Symmetrically, one gets $(I +
\tilde R(x))(I - R(x)) = I$. Hence the operator $(I + \tilde
R(x))^{-1}$ exists, and it is a linear bounded operator. \B

\medskip

Equation~(24) is called {\it the main equation} of the inverse problem.
Solving~(24) we find the vector $\psi(x)$, and consequently, the
functions $\vv_{nqi}(x)$ by formulae
$$
   \begin{array}{c}
   \vv_{n m_s 1}(x) = \psi_{n m_s 1}(x), \quad
   \vv_{n m_s 0} (x) = \vv_{n m_s 1}(x) + \xi_n \psi_{n m_s 0}(x), \\
   \vv_{n q i}(x) = \vv_{n m_s i}(x) + \xi_n \psi_{nqi}(x), \\
   n \geq 0, \, s = \overline{1, p}, \, m_s < q < m_{s+1}, \, i = 0, 1.
   \end{array}                                                                                \eqno(25)
$$
Then we construct the potential $Q(x)$ and the coefficients of the boundary conditions $h$
and $H$ via (18).
Thus, we obtain the following algorithm for the solution of
Inverse Problem~1.

\medskip

{\bf Algorithm 1.} {\it Given the data $\Lambda.$

(1) Choose $\tilde L \in A(\om)$, and calculate $\tilde \psi(x)$
and $\tilde R(x).$

(2) Find $\psi(x)$ by solving equation~(24), and calculate
$\vv_{nqi}(x).$

(3) Construct $Q(x)$, $h$ and $H$ by~(18).
}\\

{\bf 4. Sufficiency.}\\

{\it 4.1.} Let data $\{\la_{nq}, \alpha_{nq}\}_{n \geq 0, q =
\overline{1, m}} \in \mbox{Sp}$ satisfying the conditions of
Theorem~1 be given. Choose $\tilde L \in A(\om)$, construct
$\tilde \psi(x)$, $\tilde R(x)$, and consider the equation~(24).

\medskip

{\bf Lemma 9.}
{\it
For each fixed $x \in [0, \pi]$, the operator $I + \tilde R(x)$, acting from $B$ to $B$,
has a bounded inverse operator, and the main equation~(24)
has a unique solution $\psi(x) \in B$.
}

\medskip

{\bf Proof.}
It is sufficient to prove that the homogeneous equation
$$
    \beta(x) (I + \tilde R(x)) = 0,                                                           \eqno(26)
$$
where $\beta(x) = [\beta_u(x)]_{u \in V} $, $\beta_u(x)$ are $m \times m$ matrices,
has only the zero solution.
Let $\beta(x) \in B$ be a solution of (26), i.\,e.
$$
    \beta_{nqi}(x) + \sum_{(k, l, j) \in V} \beta_{klj}(x) \tilde R_{klj, nqi} (x) = 0_m.
$$
Denote
$$
    \begin{array}{c}
   \gamma_{n m_s 1}(x) = \beta_{n m_s 1}(x), \quad
   \gamma_{n m_s 0} (x) = \gamma_{n m_s 1}(x) + \xi_n \beta_{n m_s 0}(x), \\
   \gamma_{n q i}(x) = \gamma_{n m_s i}(x) + \xi_n \beta_{nqi}(x), \\
   n \geq 0, \, s = \overline{1, p}, \, m_s < q < m_{s+1}, \, i = 0, 1.
   \end{array}
$$
Then $\gamma_{nqi}(x)$ satisfy the relations
$$
    \gamma_{nqi}(x) + \sum_{k = 0}^{\iy} \sum_{l = 1}^m
    (\gamma_{kl0}(x) \tilde F_{kl0, nqi}(x) - \gamma_{kl1}(x)
    \tilde F_{kl1, nqi}(x)) = 0_m, \quad n \geq 0,                                            \eqno(27)
$$
and the following estimates are valid
$$
   \begin{array}{c}
   \|\gamma_{nqi}(x) \| \leq C(x), \quad n \geq 0, \quad q = \overline{1, m}, \\
   \|\gamma_{n {m_s} 0}(x) - \gamma_{n {m_s} 1}(x) \|,
   \|\gamma_{nqi}(x) - \gamma_{nm_si}(x) \| \leq C(x) \xi_n, \\
   s = \overline{1, p}, \, m_s < q < m_{s+1}.
   \end{array}                                                                                \eqno(28)
$$

Construct the matrix-functions $\gamma(x, \la)$, $\Gamma(x, \la)$
and $B(x, \la)$ by the formulas
$$
   \gamma(x, \la) = -\sum\limits_{k = 0}^{\iy} \sum\limits_{l = 1}^m \biggl[
   \gamma_{kl0}(x) \alpha'_{kl0} \frac{\langle \tilde \vv_{kl0}^*(x),
   \tilde \vv(x, \la) \rangle}{\la - \la_{kl0}}
   - \gamma_{kl1}(x) \alpha'_{kl1} \frac{\langle \tilde \vv_{kl1}^*(x),
   \tilde \vv(x, \la) \rangle}{\la - \la_{kl1}}
   \biggr],                                                                                   \eqno(29)
$$
$$
    \Gamma(x, \la) = -\sum\limits_{k = 0}^{\iy} \sum\limits_{l = 1}^m \biggl[
    \gamma_{kl0}(x) \alpha'_{kl0} \frac{\langle \tilde \vv_{kl0}^*(x),
    \tilde \Phi(x, \la) \rangle}{\la - \la_{kl0}}
     - \gamma_{kl1}(x) \alpha'_{kl1} \frac{\langle \tilde \vv_{kl1}^*(x),
     \tilde \Phi(x, \la) \rangle}{\la - \la_{kl1}}
    \biggr],                                                                                  \eqno(30)
$$
$$
     B(x, \la) = \gamma^*(x, \bar{\la}) \Gamma(x, \la).
$$

In view of~(12), the matrix-function $\gamma(x, \la)$ is entire in
$\la$ for each fixed $x$. The functions $\Gamma(x, \la)$ and $B(x,
\la)$ are meromorphic in  $\la$ with simple poles $\la_{nqi}$.
According to (29), $\gamma(x, \la_{nqi}) = \gamma_{nqi}(x)$. We
calculate residues of $B(x, \la)$ (for simplicity we assume that
$\{ \la_{nq0}\} \cap \{\la_{nq1} \} = \emptyset$):
$$
    \mathop{\mathrm{Res}}_{\la = \la_{nq0}} B(x, \la) =
    \gamma^*(x, \la_{nq0}) \gamma(x, \la_{nq0}) \alpha_{nq0}, \quad
    \mathop{\mathrm{Res}}_{\la = \la_{nq1}} B(x, \la) = 0_m.
$$

Consider the integral
$$
    I_N(x) = \frac{1}{2 \pi i} \int\limits_{\Gamma_N} B(x, \la) \, d \la,
$$
where $\Gamma_N = \{ \la \colon |\la| = (N + 1/2)^2\}$. Let us
show that for each fixed $x \in [0, \pi]$
$$\lim\limits_{N \to \iy} I_N(x) = 0_m.$$

Indeed, it follows from (12) and (29) that
$$
  -\gamma(x, \la) = \sum\limits_{k = 0}^{\iy}\sum\limits_{s = 1}^p \sum_{l = m_s}^{m_{s+1}-1}
  \Bigl[ \gamma_{kl0}(x) \alpha'_{kl0} \tilde D(x, \la, \la_{kl0}) -
  \gamma_{kl1}(x) \alpha'_{kl1} \tilde D(x, \la, \la_{kl1})\Bigr]
$$
$$
   =\sum_{k = 0}^{\iy}\sum_{s = 1}^p \Bigl[ (\gamma_{k m_s 0}(x) - \gamma_{k m_s 1}(x))
  \sum_{l = m_s}^{m_{s+1}-1} \alpha'_{kl0} \tilde D(x, \la, \la_{kl0}) +
  \gamma_{k m_s 1}(x) \alpha^{(s)}_k (\tilde D(x, \la, \la_{k m_s 0})
$$
$$
  - \tilde D(x, \la, \la_{k m_s 1})) +
  \gamma_{k m_s 1}(x) (\alpha^{(s)}_k - \tilde \alpha^{(s)}_k) \tilde D(x, \la, \la_{k m_s 1})
  + \gamma_{k m_s 1}(x) \sum_{l = m_s}^{m_{s+1} - 1} \sum_{j = 0}^1 \alpha'_{klj} (\tilde D(x, \la, \la_{klj})
$$
$$
  - \tilde D(x, \la, \la_{k m_s j}))
  + \sum_{l = m_s}^{m_{s+1} - 1} \sum_{j = 0}^1 (\gamma_{klj}(x) - \gamma_{k m_s j}(x)) \alpha'_{klj} \tilde D(x, \la, \la_{klj})
  \Bigr].
$$

By virtue of Lemma~6, (5) and (28), we get
$$
   \|\gamma(x, \la) \| \leq C(x) \exp(|\tau|x) \sum_{k = 0}^{\iy} \frac{\xi_k}{|\rho - k| + 1},
   \quad \mbox{Re} \, \rho \geq 0.
$$
Similarly, using (30) we obtain for sufficiently large $\rho^* > 0$:
$$
    \|\Gamma(x, \la) \| \leq \frac{C(x)}{|p|} \exp(-|\tau|x) \sum_{k = 0}^{\iy} \frac{\xi_k}{|\rho - k| + 1},
    \mbox{Re} \, \rho \geq 0, \, |\rho| \geq \rho^*, \, \rho \in G_{\de}.
$$
Then
$$
    \| B(x, \la) \| \leq \frac{C(x)} {|\rho|} \left( \sum_{k = 0}^{\iy} \frac{\xi_k}{|\rho - k| + 1} \right)^2
    \leq \frac{C(x)} {|\rho|^3}, \quad \la \in \Gamma_N.
$$
This estimate yields $\lim\limits_{N \to \iy} I_N(x) = 0_m$.

On the other hand, calculating the integral $I_N(x)$ by the residue theorem, we arrive at
$$
     \sum\limits_{k = 0}^{\iy} \sum\limits_{q = 1}^m \gamma^*_{kl0}(x) \gamma_{kl0}(x) \alpha'_{kl0} = 0_m.
$$
Since $\alpha_{kl0} = \alpha^*_{kl0} \geq 0$, we get
$$\gamma^*_{kl0}(x) \gamma_{kl0}(x) \alpha_{kl0} = 0_m,$$
$$\gamma(x, \la_{kl0}) \alpha_{kl0} = 0_m, \quad k \geq 0, \quad l = \overline{1, m}.$$

Since $\gamma(x, \la)$ is entire in $\la$, and
 $$\gamma(x, \la) = O(\exp(|\tau| x))$$
 for each fixed $x \in [0, \pi]$,
 according to condition~3 of Theorem~1, we get $\gamma(x, \la) \equiv 0_m$.
Therefore $\gamma_{nqi}(x) = 0_m$ for all $n \geq 0$, $q = \overline{1, m}$,
$i = 0, 1$, i.\,e. the homogeneous equation (26) has only the zero solution.
\B

\medskip

{\it 4.2.} Further, we provide the general strategy of the proof
of sufficiency in Theorem~1. The proofs of Lemmas~10-12 are similar
to ones described in [3, Sec.~1.6.2].

Let $\psi(x) = [\psi_u(x)]_{u \in V}$ be the solution of the main equation~(24).

\medskip

{\bf Lemma 10.}
{\it
For $n \geq 0$, $q = \overline{1, m}$, $i = 0, 1$, the following relations hold
$$
\begin{array}{c}
   \psi_{nqi}(x) \in C^1[0, \pi], \quad \| \psi^{(\nu)}_{nqi}\| \leq C(n + 1)^{\nu}, \quad \nu = 0, 1 \quad x \in [0, \pi], \\
   \|\psi_{nqi}(x) - \tilde  \psi_{nqi}(x)\| \leq C \Omega \eta_n, \quad
   \| \psi'_{nqi}(x) - \tilde \psi'_{nqi}(x) \| \leq C \Omega, \quad x \in [0, \pi],
\end{array}
$$
where
$$
    \eta_n := \left(\sum_{k = 0}^{\iy} \frac{1}{(k + 1)^2(|n - k| + 1)^2} \right).
$$
}

\medskip

Construct matrix-functions $\vv_{nqi}(x)$ by formulae~(25). By
virtue of Lemma~10, we have
$$
    \begin{array}{c}
    \| \vv_{nqi}^{(\nu)}(x)\| \leq C(n + 1)^{\nu}, \quad \nu = 0, 1, \\
    \| \vv_{nqi}(x) - \tilde \vv_{nqi}(x)\| \leq C \Omega \eta_n, \quad
    \| \vv'_{nqi}(x) - \tilde \vv'_{nqi}(x)\| \leq C \Omega, \quad q = \overline{1, m}, \\
    \| \vv_{n m_s 0}(x) - \vv_{n m_s 1}(x) \|, \, \| \vv_{n q i}(x) - \vv_{n m_s i}(x)\| \leq C \xi_n, \quad
     s = \overline{1, p}, \, m_s < q < m_{s+1}.
     \end{array}                                                                              \eqno(31)
$$

Further, we construct the matrix-functions $\vv(x, \la)$ and
$\Phi(x, \la)$ by the formulas
$$
    \vv(x, \la) = \tilde \vv(x, \la) - \sum_{(k, l, j) \in V} (-1)^j
    \vv_{klj}(x) \alpha'_{klj} \frac{\langle \tilde \vv^*_{klj}(x),  \tilde \vv(x, \la) \rangle}{\la - \la_{klj}},
$$
$$
    \Phi(x, \la) = \tilde \Phi(x, \la) - \sum_{(k, l, j) \in V} (-1)^j
    \vv_{klj}(x) \alpha'_{klj} \frac{\langle \tilde \vv^*_{klj}(x),  \tilde \Phi(x, \la) \rangle}{\la - \la_{klj}},
$$
and the boundary value problem $L(Q(x), h, H)$ via (18). Clearly,
$\vv(x, \la_{nqi}) = \vv_{nqi}(x)$.

Using estimates~(31) one can show that the components of $\ee_0(x)$ are
absolutely continuous and the components of
$\ee(x)$ belong to $L_2(0, \pi)$.
Consequently, we get

\medskip

{\bf Lemma 11.}
{\it
$Q_{jk}(x) \in L_2(0, \pi)$, $j, k = \overline{1, m}$.
}

\medskip

{\bf Lemma 12.}
{\it
The following relations hold
$$
   \ell \vv_{nqi}(x) = \la_{nqi} \vv_{nqi}(x), \quad
   \ell \vv(x, \la) = \la \vv(x, \la), \quad \ell \Phi(x, \la) = \la \Phi(x, \la),
$$
$$
   \vv(0, \la) = I_m, \quad \vv'(0, \la) = h, \quad U(\Phi) = I_m, \quad V(\Phi) = 0_m.
$$
}

\medskip

In order to finish the proof of Theorem~1 it remains to show that
the given data $\{\la_{nq}, \alpha_{nq} \}$ coincide with the
spectral data of the constructed boundary value problem
 $L(Q, h, H)$.
In view of Lemma~12, the matrix-function $\Phi(x, \la)$ is the
Weyl solution of $L$. Let us get the representation for the Weyl
matrix:
$$
    M(\la) = \Phi(0, \la) = \tilde M(\la) -  \sum_{(k, l, j) \in V} \vv_{klj}(0) \alpha'_{klj}
    \frac{\langle \tilde \vv^*_{klj}(x), \tilde \Phi(x, \la) \rangle_{x = 0}}{ \la - \la_{klj}}
    \tilde M(\la)
$$
$$
    + \sum_{k = 0}^{\iy} \sum_{l = 1}^m \left(\frac{\alpha'_{kl0}}{\la - \la_{kl1}} - \frac{\alpha'_{kl1}}{\la - \la_{kl1}} \right).
$$
Using the equality (see [4])
$$
    \tilde M(\la) = \sum_{k = 0}^{\iy} \sum_{l = 1}^m \frac{\alpha'_{kl1}}{\la - \la_{kl1}},
$$
we arrive at
$$
     M(\la) = \sum_{k = 0}^{\iy} \sum_{l = 1}^m \frac{\alpha'_{kl0}}{\la - \la_{kl0}}.
$$
Consequently, $\{ \la_{kl0} \}$ are simple poles of the Weyl
matrix $M(\la)$, and $\{ \alpha_{kl0}\}$ are residues at the
poles. Note that the multiplicities of the eigenvalues coincide
with the numbers of equal values among $\{ \la_{kl0}\}$, because
they both coincide with the ranks of $\{ \alpha_{kl0}\}$.
Theorem~1 is proved. \B

\medskip

{\bf Acknowledgment.}  This research was supported in part by Grants
10-01-00099 and 10-01-92001-NSC of Russian Foundation for Basic Research
and Taiwan National Science Council and by the Moebius Contest Foundation for Young
Scientists.

\begin{center}
{\bf REFERENCES}
\end{center}
\begin{enumerate}
\item[{[1]}] Marchenko V.A. Sturm-Liouville Operators and their Applications,
     Naukova Dumka, Kiev, 1977 (Russian); English transl., Birkhauser, 1986.
\item[{[2]}] Levitan B.M. Inverse Sturm-Liouville Problems, Nauka, Moscow,
     1984 (Russian); English transl., VNU Sci.Press, Utrecht, 1987.
\item[{[3]}] Freiling G. and Yurko V.A. Inverse Sturm-Liouville Problems
     and their Applications. NOVA Science Publishers, New York, 2001.
\item[{[4]}] Yurko V.A. Inverse problems for matrix Sturm-Liouville operators,
    Russian J. Math. Phys. 13, no.1 (2006), 111-118.
\item[{[5]}] Carlson R. An inverse problem for the matrix Schr\"{o}dinger equation,
    J. Math. Anal. Appl. 267 (2002), 564-575.
\item[{[6]}] Malamud M.M. Uniqueness of the matrix Sturm-Liouville equation
    given a part of the monodromy matrix, and Borg type results.
    Sturm-Liouville Theory, Birkh\"{a}user, Basel, 2005, 237-270.
\item[{[7]}] Yurko V.A. Inverse problems for the matrix Sturm-Liouville equation on a finite interval,
    Inverse Problems, 22 (2006), 1139-1149.
\item[{[8]}] Chelkak D., Korotyaev E. Weyl-Titchmarsh functions of vector-valued
    Sturm-Liouville operators on the unit interval,
    J. Func. Anal. 257 (2009), 1546-1588.
\item[{[9]}] Mykytyuk Ya.V., Trush N.S. Inverse spectral problems for Sturm-Liouville
    operators with matrix-valued potentials,
    Inverse Problems, 26 (2010), 015009.
\item[{[10]}] Naimark M.A. Linear Differential Operators, 2nd ed., Nauka,
     Moscow, 1969; English transl. of 1st ed., Parts I,II, Ungar,
     New York, 1967, 1968.
\end{enumerate}

\vspace{1cm}

Natalia Bondarenko

Department of Mathematics

Saratov State University

Astrakhanskaya 83, Saratov 410026, Russia

bondarenkonp@info.sgu.ru

\end{document}